\documentclass[11pt]{article}
\newcommand{\Min}{\operatorname{Min}}
\newcommand{\MIN}{\operatorname{MIN}}
\newcommand{\COL}{\operatorname{COL}}
\newcommand{\Max}{\operatorname{Max}}
\newcommand{\MAX}{\operatorname{MAX}}
\newcommand{\ROW}{\operatorname{ROW}}
\newcommand{\SB}{{\cal SB}}
\renewcommand{\sb}{{\mathrm{sb}}}
\newcommand{\CB}{{\cal CB}}
\newcommand{\cb}{{\mathrm{cb}}}
\newcommand{\PM}{\operatorname{PM}}
\newcommand{\fr}[1]{\stackrel{\frown}{#1}}
\newcommand{\KG}{K_{\mathbb G}}
\renewcommand{\baselinestretch}{1.2}
\newcommand{\dated}{\mbox{} \hfill {\small [{\tt \today}]}} 
\newtheorem{theorem}{Theorem}[section]
\newtheorem{lemma}[theorem]{Lemma}
\newtheorem{corollary}[theorem]{Corollary}
\newtheorem{proposition}[theorem]{Proposition}
\newtheorem{df}[theorem]{Definition}
\newenvironment{definition}{\begin{df} \rm}{\end{df}}

 \usepackage{amsmath,amssymb,amscd}
\newcommand{\pf}[1]{\trivlist \item[\hskip\labelsep\it #1\ ]}
\newcommand{\varpf}[1]{\trivlist \item[\hskip\labelsep\sc #1:]}
\newcommand{\qedbox}{$\rlap{$\sqcap$}\sqcup$}
\newcommand{\qed}{\qquad \qedbox \endtrivlist}
\newcommand{\varqed}{\hfill \rule{0.6em}{0.6em} \endtrivlist}
\newenvironment{proof}{\pf{Proof}}{\qed}

\newenvironment{remark}{\pf{Remark}}{\endtrivlist}
\newenvironment{remarks}{\pf{Remarks} 
   \begin{enumerate}}{\end{enumerate} \endtrivlist}

\newenvironment{examples}{\pf{Examples} 
   \begin{enumerate}}{\end{enumerate} \endtrivlist}
\newenvironment{items}{
  \begin{enumerate} 
                    
  }{\end{enumerate}}
\newenvironment{alphitems}{
  \begin{enumerate} 
                    
  }{\end{enumerate}}
\newenvironment{keywords}{\noindent\small {\it Keywords\/}:}{\vskip 4pt}
\newenvironment{classification}{\noindent\small 2000 {\it Mathematics Subject
Classification\/}:}{\vskip 12pt}

\renewcommand{\iff}{\quad\Longleftrightarrow\quad}

\newcommand{\comps}{{\mathbb C}}

\newcommand{\posints}{{\mathbb N}}

\newcommand{\tensor}{\otimes}

\newcommand{\Tensor}{\hat{\otimes}}

\newcommand{\wTensor}{\check{\otimes}}

\newcommand{\cstar}{{C^\ast}}

\newcommand{\id}{{\mathrm{id}}}

\newcommand{\diag}{{\operatorname{diag}}}

\newcommand{\A}{{\mathfrak A}}
\newcommand{\B}{{\mathfrak B}}
\newcommand{\Hilbert}{{\mathfrak H}}
\newcommand{\M}{{\mathfrak M}}

\newcommand{\VN}{\operatorname{VN}}

\title{Operator space structure and amenability \\ for Fig\`a-Talamanca--Herz algebras}
\author{{\it Anselm Lambert} \and {\it Matthias Neufang}\footnote{Part of the research for this paper was done while the author held a postdoctoral fellowship of the Pacific Institute for the Mathematical Sciences
(PIMS) at the University of Alberta; the support by PIMS is acknowledged with thanks.} \and {\it Volker Runde}\thanks{Research supported by NSERC under grant no.\ 227043-00.}}
\date{}
\begin{document}
\maketitle
\begin{abstract}
Column and row operator spaces --- which we denote by $\COL$ and $\ROW$, respectively --- over arbitrary Banach spaces were introduced by the first-named author; for Hilbert spaces, these definitions coincide with the usual ones. 
Given a locally compact group $G$ and $p,p' \in (1,\infty)$ with $\frac{1}{p} + \frac{1}{p'} = 1$, we use the operator space structure on $\CB(\COL(L^{p'}(G)))$ to equip the Fig\`a-Talamanca--Herz algebra $A_p(G)$ 
with an operator space structure, turning it into a quantized Banach algebra. Moreover, we show that, for $p \leq q \leq 2$ or $2 \leq q \leq p$ and amenable $G$, the canonical inclusion $A_q(G) \subset A_p(G)$ is 
completely bounded (with $\cb$-norm at most $\KG^2$, where $\KG$ is Grothendieck's constant). As an application, we show that $G$ is amenable if and only if $A_p(G)$ is operator amenable for all --- and equivalently 
for one --- $p \in (1,\infty)$; this extends a theorem by Z.-J.\ Ruan. 
\end{abstract}
\section*{Introduction}
The Fourier algebra $A(G)$ of a locally compact group $G$ was introduced by P.\ Eymard in \cite{Eym1}. If $G$ is abelian with dual group $\Gamma$, then the Fourier transform induces an isometric
isomorphism of $A(G)$ and $L^1(\Gamma)$. Although the Fourier algebra is an invariant for $G$ --- like $L^1(G)$ ---, its Banach algebraic amenability does not correspond well to the amenability of
$G$ --- very much unlike $L^1(G)$: The group $G$ is amenable if and only if $L^1(G)$ is amenable as a Banach algebra (\cite{Joh1}), but there are compact groups, among them $\operatorname{SO}(3)$, for which $A(G)$ 
fails to be even weakly amenable (\cite{Joh2}). In fact, the only locally compact groups $G$ for which $A(G)$ is an amenable Banach algebra in the sense of \cite{Joh1} are those with an abelian subgroup of finite index 
(\cite{Run3}).
\par
Since $A(G)$ is the predual of the group von Neumann algebra $\VN(G)$, it is an operator space in a natural manner. In \cite{Rua1}, Z.-J.\ Ruan introduced a variant of amenability for ``quantized'' Banach algebras --- 
called operator amenability --- that reflects the operator space structures of those algebras. He showed that a locally compact group $G$ is amenable if and only if $A(G)$ is operator amenable (\cite[Theorem 3.6]{Rua1}). 
Further results by O.\ Yu.\ Aristov (\cite{Ari}), P.\ J.\ Wood (\cite{Woo}), and N.\ Spronk (\cite{Spr}) lend additional support to the belief that homological properties of $A(G)$, such as amenability,
biprojectivity or weak amenability, correspond to properties of $G$ much more naturally if the operator space structure is taken into account. Even if one is only interested in $A(G)$ as a Banach algebra, considering the
canonical operator space structure can be valuable: Although the main result of \cite{Run3} is purely classical in its statement, its proof is operator space theoretic. (Further examples of classical results obtained with the help of
operator space methods can be found in \cite{FKLS}.)
\par
In \cite{Her1}, C.\ Herz introduced, for $p \in (1,\infty)$, an $L^p$-analog of the Fourier algebra, denoted by $A_p(G)$. These algebras are called
Fig\`a-Talamanca--Herz algebras. It was asked by the third-named author if there was an analog of Ruan's theorem for arbitrary Fig\`a-Talamanca--Herz algebras (\cite[Problem 34]{LoA}). 
The first obstacle towards a solution of this problem is that there is --- at first glance --- no natural operator space structure for $A_p(G)$ if $p \in (1,\infty) \setminus \{ 2 \}$. 
\par
In \cite{Run}, the same author used the operator space structure introduced by G.\ Pisier on the $L^p$-spaces via complex interpolation (\cite{Pis}) to define
operator space analogs $OA_p(G)$ of the classical Fig\`a-Talamanca--Herz algebras for all $p \in (1,\infty)$. These operator Fig\`a-Talamanca--Herz algebras display properties similar to those of their classical
counterparts (see, for example, \cite[Theorem 4.10]{Run}). Nevertheless, the construction from \cite{Run} is unsatisfactory for two reasons:
\begin{enumerate}
\item Although we have a contractive inclusion $A_p(G) \subset OA_p(G)$ for all $p \in (1,\infty)$, the two algebras $A_p(G)$ and $OA_p(G)$ can fail to be identical as Banach spaces for $p \neq 2$.
\item Even though $A(G) = OA_2(G)$ as Banach spaces, they need not have the same operator space structure.
\end{enumerate}
The reason why $A(G)$ and $OA_2(G)$ need not coincide as operator spaces is that the operator space structure of $A(G)$ stems from $\VN(G)$ and thus from the column Hilbert space over $L^2(G)$, whereas
$OA_2(G)$ is defined by means of the operator Hilbert space over $L^2(G)$. In order to get a satisfactory operator space structure for general Fig\`a-Talamanca--Herz algebras, one should therefore strive for an 
extension of the notion of column space from Hilbert spaces to arbitrary $L^p$-spaces. 
\par
In his doctoral dissertation \cite{Lam}, the first-named author introduced such a notion; in fact, he defined column and row spaces over arbitrary Banach spaces. We first outline his construction and then use it  to equip general 
Fig\`a-Talamanca--Herz algebras with a canonical operator space structure, turning them into ``quantized'' Banach algebras. As an application, we extend \cite[Theorem 3.6]{Rua1} to arbitrary Fig\`a-Talamanca--Herz algebras.
\subsection*{Acknowledgments}
The authors would like to thank Gerd Wittstock for several valuable discussions on the topic of this paper. The first-named author wishes to thank him especially for
his guidance and supervision that ultimately resulted in \cite{Lam}. 
\section{Preliminaries}
\subsection{Fig\`a-Talamanca--Herz algebras}
Let $G$ be a locally compact group. For any function $f \!: G \to \comps$, we define $\check{f} \!: G \to \comps$ by letting $\check{f}(x) := f(x^{-1})$ for $x
\in G$. Let $p \in (1,\infty)$, and let $p' \in (1,\infty)$ be dual to $p$, i.e.\ $\frac{1}{p} + \frac{1}{p'} = 1$. The {\it Fig\`a-Talamanca--Herz algebra\/} $A_p(G)$ consists of those functions $f \!: G \to \comps$
such that there are sequences $( \xi_n )_{n=1}^\infty$ in $L^p(G)$ and $( \eta_n )_{n=1}^\infty$ in $L^{p'}(G)$ such that
\begin{equation} \label{Apeq1}
  \sum_{n=1}^\infty \| \xi_n \|_{L^p(G)} \| \eta_n \|_{L^{p'}(G)} < \infty
\end{equation}
and
\begin{equation} \label{Apeq2}
  f = \sum_{n=1}^\infty \xi_n \ast \check{\eta}_n.
\end{equation}
The norm on $A_p(G)$ is defined as the infimum over all sums (\ref{Apeq1}) such that (\ref{Apeq2}) holds. It is clear that $A_p(G)$ is a Banach space that embeds contractively into ${\cal C}_0(G)$. 
It was shown by C.\ Herz (\cite{Her1}) that $A_p(G)$ is closed under pointwise multiplication and, in fact, a Banach algebra. The case where $p = q = 2$ had previously been studied by P.\ Eymard (\cite{Eym1}); in this case 
$A(G) := A_2(G)$ is called the {\it Fourier algebra\/} of $G$.
\par
Let $\lambda_{p'} \!: G \to {\cal B}(L^{p'}(G))$ be the regular left representation of $G$ on $L^{p'}(G)$. Via integration, $\lambda_{p'}$ extends to a representation of $L^1(G)$ on $L^{p'}(G)$. 
The algebra of {\it $p'$-pseudomeasures\/} $\PM_{p'}(G)$ is defined as the $w^\ast$-closure of $\lambda_{p'}(L^1(G))$ in ${\cal B}(L^{p'}(G))$.
There is a canonical duality $\PM_{p'}(G) \cong A_p(G)^\ast$ via
\[
  \langle \xi \ast \check{\eta}, T \rangle := \langle T\eta, \xi \rangle \qquad (\xi \in L^{p'}(G), \, \eta \in L^p(G), \, T \in \PM_{p'}(G)). 
\]
If $p =2$, then $\VN(G) := \PM_2(G)$ is a von Neumann algebra, the {\it group von Neumann algebra\/} of $G$. 
\par
For more information, see \cite{Eym1}, \cite{Eym2}, \cite{Her1}, \cite{Her2}, and \cite{Pie}.
\subsection{Operator spaces and quantized Banach algebras}
There is now a booklength monograph available on the subject of operator spaces (\cite{ER}), and a second one will soon appear (\cite{PisBk}); furthermore, a very accessible survey article is available on the internet 
(\cite{Wit}). We therefore refrain from repeating here the basic definitions of operator
space theory. In our choice of notation, we mostly follow \cite{ER}. In particular, the projective and injective tensor product of Banach spaces are denoted by $\tensor^\gamma$ and $\tensor^\lambda$, respectively, whereas
$\Tensor$ and $\wTensor$ stand for the projective and injective tensor product of operator spaces. 
\par
We deviate from \cite{ER} with our notation in two points:
\begin{enumerate}
\item The minimal and maximal operator space over a given Banach space $E$ is denoted by $\MIN(E)$ and $\MAX(E)$, respectively.
\item The column and row space over a Hilbert space $\Hilbert$ is denoted by $\COL(\Hilbert)$ and $\ROW(\Hilbert)$, respectively.
\end{enumerate}
\begin{definition}
A {\it quantized Banach algebra\/} is an algebra which is also an operator space such that multiplication is completely bounded.
\end{definition}
\begin{remark}
We do not require multiplication to be completely contractive (such algebras are called {\it completely contractive Banach algebras\/}; see \cite{Rua1}). In our choice of terminology, we follow \cite{Ari}.
\end{remark}
\begin{examples}
\item For any Banach algebra $\A$ --- not necessarily with contractive multiplication --- the maximal operator space $\MAX(\A)$ is a quantized Banach algebra.
\item If $\mathfrak H$ is a Hilbert space, then any closed subalgebra of ${\cal B}({\mathfrak H})$ is completely contractive.
\item We denote the $W^\ast$-tensor product by $\bar{\tensor}$. A {\it Hopf--von Neumann algebra\/} is a pair $(\M, \nabla)$, where $\M$ is a von Neumann algebra, and 
$\nabla$ is a {\it co-multiplication\/}: a unital, $w^\ast$-continuous, and injective $^\ast$-homomorphism $\nabla \!: \M \to \M \bar{\tensor} \M$ which is co-associative, i.e.\ the diagram
\[
  \begin{CD}
  \M  @>\nabla>> \M \bar{\tensor} \M \\
  @V{\nabla}VV                                          @VV{\nabla \tensor \id_\M}V                     \\
  \M \bar{\tensor} \M @>>{\id_\M \tensor \nabla}> \M \bar{\tensor} \M \bar{\tensor} \M
  \end{CD}    
\]
commutes. Let $\M_\ast$ denote the unique predual of $\M$. By \cite[Theorem 7.2.4]{ER}, we have $\M \bar{\tensor} \M \cong (\M_\ast \Tensor \M_\ast)^\ast$. Thus $\nabla$ induces a
complete contraction $\nabla_\ast \!: \M_\ast \Tensor \M_\ast \to \M_\ast$ turning $\M_\ast$ into a completely contractive Banach algebra. 
\item Let $G$ be a locally compact group. Then the representation
\[
  G \mapsto {\cal B}(L^2(G \times G)), \quad x \mapsto \lambda_2(x) \tensor \lambda_2(x)
\]
induces a co-multiplication $\nabla \!: \VN(G) \to \VN(G \times G) \cong \VN(G) \bar{\tensor} \VN(G)$. Hence, $A(G)$ with its canonical operator space structure is a completely contractive Banach algebra.
\end{examples}
\section{Operator sequence spaces}
In \cite{Math}, B.\ Mathes characterized the column operator space $\COL(\Hilbert)$ over a Hilbert space $\Hilbert$ through the isometries
\[
  M_{n,1}(\COL(\Hilbert)) = M_{n,1}(\MAX(\Hilbert)) \qquad\text{and}\qquad M_{1,n}(\COL(\Hilbert)) = M_{1,n}(\MIN(\Hilbert))
\]
for all $n \in \posints$, i.e.\ $\COL(\Hilbert)$ is maximal on the
columns and minimal on the rows. A similar characterization holds for
$\ROW(\Hilbert)$.
\par
In order to define column and row operator spaces over {\it
  arbitrary\/} Banach spaces in the next section, we first introduce
and discuss an axiomatic characterization of the columns of operator
spaces: the operator sequence spaces. Those spaces were introduced by the first-named author in his doctoral dissertation (\cite[Definition 1.1.1]{Lam}). They form a category somewhere in between Banach and operator spaces.
A full account of the theory of operator sequence spaces will be published elsewhere. 
\par
In this section, we content ourselves with presenting the concepts and
results we need for the remainder of the paper: to define column
and row operator spaces over arbitrary Banach spaces and to use those operator space
structures to turn the Fig\`a-Talamanca--Herz algebras into quantized
Banach algebras. We are somewhat sketchy with our proofs --- especially if they consist mainly of routine calculations or straightforward adaptations 
of proofs of the corresponding Banach or operator space results.
\begin{definition}
A {\it sequential norm\/} over a Banach space $E$ is a sequence $( \| \cdot \|_{\fr{n}} )_{n=1}^\infty$ such that $\| \cdot \|_{\fr{1}}$ is the given norm on $E$ and $\| \cdot \|_{\fr{n}}$ is, for
each $n \in \posints$, a norm on $E^n$ such that
\begin{eqnarray*}
  \left\| \left[ \begin{array}{c} x \\ 0 \end{array} \right] \right\|_{\fr{m+n}} & = & \| x \|_{\fr{m}} \qquad (m,n \in \posints, \, x \in E^m), \\
  \left\| \left[ \begin{array}{c} x \\ y \end{array} \right] \right\|_{\fr{m+n}}^2 & \leq & \| x \|_{\fr{m}}^2 + \| y \|_{\fr{n}}^2 \qquad (m,n \in \posints, \, x \in E^m, \, y \in F^n)
\end{eqnarray*} 
and
\[
  \| \alpha x \|_{\fr{m}} \leq \| \alpha \| \| x \|_{\fr{n}} \qquad (m,n \in \posints, \, x \in E^n, \, \alpha \in M_{m,n}).
\]
For $n \in \posints$, we write $E^{\fr{n}}$ to denote $E^n$ equipped with $\| \cdot \|_{\fr{n}}$. The space $E$ together with the sequential norm $( \| \cdot \|_{\fr{n}} )_{n=1}^\infty$ is called
an {\it operator sequence space\/} (over $E$).
\end{definition}
\begin{examples}
\item Letting $\comps^{\fr{n}} := \ell^2_n$ for $n \in \posints$, we define {\it the\/} (unique) operator sequence space over $\comps$.
\item Let $E$ be any Banach space. The {\it minimal operator sequence space\/} $\min(E)$ over $E$ is defined by letting $\min(E)^{\fr{n}} := {\cal B}(\ell^2_n,E)$ for $n \in \posints$; the adjective minimal will be justified below.
\item Let $E$ be any Banach space. The {\it maximal operator sequence space\/} $\max(E)$ over $E$ is defined as follows: For $n \in \posints$ and $x \in E^n$, define
\[
  \| x \|_{\fr{n}} := \inf \{ \| \alpha \| \| y \|_{\ell^2_m(E)} : m \in \posints, \, \alpha \in M_{n,m}, \, y \in E^m, \, x = \alpha y \}.
\]
As in the case of $\min(E)$, the use of the adjective maximal will soon become clear.
\item Let $E$ be an operator sequence space, and let $m \in \posints$. We define $(E^{\fr{m}})^{\fr{n}} := E^{\fr{mn}}$ for $n \in \posints$. This turns $E^{\fr{m}}$ into an operator sequence space.
\item Let $E$ be an operator space. Define an operator sequence space $C(E)$ over $M_1(E)$ by letting $C(E)^{\fr{n}} := M_{n,1}(E)$ for $n \in \posints$. It will become apparent in the next section that
{\it every\/} operator sequence space occurs in that fashion.
\end{examples}
\par
Having introduced the objects of the category of operator sequence spaces, we now turn to defining its morphisms.
\par
Given two operator sequence spaces $E$ and $F$, a linear map $T \!: E \to F$, and $n \in \posints$, let 
\[
  T^{\fr{n}} \!: E^{\fr{n}} \to F^{\fr{n}}, \quad [ x_j ]_{j=1, \ldots, n} \mapsto [ Tx_j ]_{j=1, \ldots, n}
\]
denote the $n$-th amplification of $T$.
\begin{definition} \label{sbdef}
Let $E$ and $F$ be operator sequence spaces. Then $T \in {\cal B}\left(E^{\fr{1}},F^{\fr{1}}\right)$ is called {\it sequentially bounded\/} if
\[
  \| T \|_\sb := \sup_{n \in \posints} \left\| T^{\fr{n}} \right\|_{{\cal B}\left(E^{\fr{n}}, F^{\fr{n}}\right)} < \infty.
\]
If $\| T \|_\sb \leq 1$, we call $T$ a {\it sequential contraction\/}, and if $T^{\fr{n}}$ is an isometry for each $n \in \posints$, we call $T$ a {\it sequential isometry\/}.
The collection of all sequentially bounded maps from $E$ to $F$ is denoted by $\SB(E,F)$.
\end{definition}
\begin{remarks}
\item It is straightforward that $\| \cdot \|_\sb$ turns $\SB(E,F)$ into a Banach space.
\item We write $\SB(E)$ instead of $\SB(E,E)$.
\end{remarks}
\begin{examples}
\item Let $E$ be an operator sequence space, and let $F$ be a Banach space. Then 
\begin{equation} \label{minsb}
  \SB(E,\min(F)) = {\cal B}\left(E^{\fr{1}},F\right)
\end{equation} 
holds isometrically. (This justifies the name minimal operator space.)
\item Let $E$ be a Banach space, and let $F$ be an operator sequence space. Then 
\begin{equation} \label{maxsb}
  \SB(\max(E),F) = {\cal B}\left(E,F^{\fr{1}}\right)
\end{equation}
holds isometrically. (This justifies the adjective maximal.)
\item Let $E$ be an infinite-dimensional Banach space. Then $\id_E \!: \MIN(E) \to \MAX(E)$ is not completely bounded (\cite[Theorem 2.12]{Pau1}). Interestingly, the situation for operator sequence spaces and sequentially 
bounded maps is different: For example, if $\A$ is a $\cstar$-algebra, then $\A$ is subhomogeneous if and only if $\id_\A \in \SB(\min(\A), \max(\A))$ (\cite[Satz 2.2.25]{Lam}).
\end{examples}
\par
Next, we sketch the duality theory for operator sequence spaces. We first introduce a canonical operator sequence space structure over the Banach space of all sequentially bounded maps between
operator sequence spaces:
\begin{proposition}
Let $E$ and $F$ be operator sequence spaces. Letting
\[
  \SB(E,F)^{\fr{n}} := \SB\left(E,F^{\fr{n}}\right) \qquad (n \in \posints)
\]
defines an operator sequence space over $\SB(E,F)$.
\end{proposition}
\par
We skip the proof which parallels the one of the corresponding result for operator spaces.
\par
We require an analog of Smith's lemma (\cite[Proposition 2.2.2]{ER}) for operator sequence spaces:
\begin{lemma}[Smith's lemma for operator sequence spaces] \label{smith}
Let $E$ and $F$ be operator sequence spaces such that $\dim F = m < \infty$. Then ${\cal B}\left(E^{\fr{1}},F^{\fr{1}}\right) = \SB(E,F)$ holds with $\| T \|_\sb = \left\| T^{\fr{m}} \right\|$ for each $T \in {\cal B}(E,F)$.
\end{lemma}
\begin{proof}
Let $T \in  {\cal B}(E,F)$, and let $n \in \posints$. Let $x = [x_j ]_{j=1, \ldots, n} \in E^{\fr{n}}$, and let $y_1, \ldots, y_m$ be a basis for $F$. 
Then there are $\alpha_{j,k} \in \comps$ for $j =1, \ldots, n$ and $k =1, \ldots, m$ such that
\[
  T x_j = \sum_{k=1}^m \alpha_{j,k} y_k.
\]
Let $y := [ y_k ]_{k=1, \ldots, m} \in F^{\fr{m}}$ and $\alpha := [ \alpha_{j,k} ]_{{j=1, \ldots, n} \atop {k=1, \ldots, m}} \in M_{n,m}$. It follows that
\[
  T^{\fr{n}} x = \alpha y.
\]
Let $v \in M_{n,m}$ be a partial isometry such that $\alpha = v |\alpha|$, where $| \alpha | := (\alpha^\ast \alpha)^\frac{1}{2}$. It follows that
\begin{equation} \label{smith1}
  \| \alpha y \|_{F^{\fr{n}}} = \| v |\alpha| y \|_{F^{\fr{n}}} \leq \| |\alpha| y \|_{F^{\fr{m}}}
\end{equation}
and
\begin{equation} \label{smith2}
  \| |\alpha| y \|_{F^{\fr{m}}} = \| v^\ast \alpha y \|_{F^{\fr{m}}} \leq \| \alpha y \|_{F^{\fr{n}}}
\end{equation}
so that in (\ref{smith1}) and (\ref{smith2}) equality holds. Consequently, we have
\begin{eqnarray*}
  \left\| T^{\fr{n}} x \right\|_{F^{\fr{n}}} & = & \| \alpha y \|_{F^{\fr{n}}} \\
  & = & \| v^\ast \alpha y \|_{F^{\fr{m}}} \\
  & = & \left\| v^\ast T^{\fr{n}} x \right\|_{F^{\fr{m}}} \\
  & = & \left\| T^{\fr{m}}( v^\ast x) \right\|_{F^{\fr{m}}} \\
  & \leq & \left\| T^{\fr{m}} \right\| \| v^\ast x \|_{F^{\fr{m}}} \\
  & \leq & \left\| T^{\fr{m}} \right\| \| x \|_{F^{\fr{n}}},
\end{eqnarray*}
so that $\left\| T^{\fr{n}} \right\| \leq \left\| T^{\fr{m}} \right\|$. Since $n \in \posints$ was arbitrary, this yields $\| T \|_\sb \leq \left\| T^{\fr{m}} \right\|$. The converse inequality is trivial.
\end{proof}
\begin{remark}
In view of \cite[Proposition 2.2.2]{ER}, one might have suspected that Lemma \ref{smith} holds true for $F = \ell^2_m$ only. The fact that there is a stronger version of Smith's lemma for operator sequence spaces
than for operator spaces has interesting consequences. For example, the principle of local reflexivity, which is a cornerstone of the local theory of Banach spaces (\cite[6.6]{DF}), but fails to have 
an analog for operator spaces (\cite[Corollary 14.3.8]{ER}), still works in the category of operator sequence spaces (\cite[Satz 1.3.26]{Lam})
\end{remark}
\begin{corollary} \label{smithcor}
Let $E$ be an operator sequence space. Then $E^\ast = \SB(E,\comps)$ holds isometrically.
\end{corollary}
\par
With Corollary \ref{smithcor} at hand, we can now equip the (Banach space) dual of an operator sequence space with a canonical operator sequence space structure.
\par
Taking the adjoint of a sequentially bounded operator yields again a sequentially bounded operator. But more is true:
\begin{theorem} \label{adjoints}
Let $E$ and $F$ be operator sequence spaces. Then
\[
  \left\| (T^\ast)^{\fr{n}} \right\| = \left\| T^{\fr{n}} \right\| \qquad ( n \in \posints, \, T \in \SB(E,F))
\]
holds. Moreover, 
\begin{equation} \label{adjoint}
  \SB(E,F) \to \SB(F^\ast,E^\ast), \quad T \mapsto T^\ast
\end{equation}
is a sequential isometry.
\end{theorem}
\begin{proof}
The first part of the theorem is \cite[Satz 1.3.14]{Lam} and has a proof similar to its operator space analog \cite[Proposition 3.2.2]{ER}. In particular, (\ref{adjoint}) is an isometry.
\par
To see that (\ref{adjoint}) is in fact a sequential isometry, fix $n \in \posints$ and note that, we have a (sequential) isometric canonical isomorphism 
\begin{equation} \label{ident}
  \SB\left(\left(F^{\fr{n}}\right)^\ast,E^\ast\right) \cong \SB\left(F^\ast,(E^\ast)^{\fr{n}}\right)
\end{equation}
by \cite[Satz 1.3.10 and Satz 1.3.12]{Lam}. Hence, we have the following canonical isometries:
\begin{eqnarray*}
  \SB(E,F)^{\fr{n}} & = & \SB\left(E,F^{\fr{n}}\right) \\
  & \hookrightarrow & \SB\left(\left(F^{\fr{n}}\right)^\ast,E^\ast\right), \qquad\text{by \cite[Satz 1.3.14]{Lam}}, \\
  & \cong & \SB\left(F^\ast,(E^\ast)^{\fr{n}}\right), \qquad\text{by (\ref{ident})}, \\
  & = & \SB(F^\ast,E^\ast)^{\fr{n}}.
\end{eqnarray*}
This completes the proof.
\end{proof}
\par
We conclude this section with the analog of the $\MIN$-$\MAX$ duality (\cite[Satz 2.1.11]{Lam}):
\begin{theorem}[$\min$-$\max$ duality] \label{minmax}
For any Banach space $E$, we have the sequentially isometric isomorphisms
\[
  \min(E)^\ast = \max(E^\ast) \qquad\text{and}\qquad \max(E)^\ast = \min(E^\ast)
\]
\end{theorem}
\begin{proof}
Since the compatibility of biduals for operator spaces (\cite[Theorem 2.5]{Ble})
has an analog in the category of operator sequence spaces (\cite[Satz
  1.3.19]{Lam}) with an almost identical proof, the $\min$-$\max$
duality can be proven by more or less following the proof of the
$\MIN$-$\MAX$ duality in \cite{Ble}.
\end{proof}
\section{Column and row spaces over arbitrary Banach spaces}
With the preparations made in the previous section, we can now define
column and row spaces over arbitrary Banach spaces. As before, the material is from \cite{Lam} and will appear elsewhere in fuller detail.
\begin{definition} \label{Mindef}
Let $E$ be an operator sequence space. Then the {\it minimal operator space\/} $\Min(E)$ over $E$ is defined by letting $M_n(\Min(E)) := {\cal B}\left( \ell^2_n, E^{\fr{n}} \right)$.
\end{definition}
\begin{remarks}
\item By \cite[Satz 4.1.2]{Lam}, $\Min(E)$ is an operator space for any operator sequence space $E$.
\item For any operator space $E$, and for any operator sequence space $F$, 
\[
  \CB(E,\Min(F)) = \SB(C(E),F)
\]
holds isometrically (\cite[Satz 4.1.6]{Lam}).
\item Let $F$ be a Banach space. Then the previous remark and the isometric identity (\ref{minsb}) combined yield that
\[
  \CB(E,\Min(\min(F))) = {\cal B}(M_1(E),F)
\]
holds isometrically for each operator space $E$, so that we have $\Min(\min(F)) = \MIN(F)$.
\end{remarks}
\par
The following definition generalizes V.\ I.\ Paulsen's formula for the maximal operator space norm over a Banach space (\cite[Theorem 2.1]{Pau2}):
\begin{definition} \label{Maxdef}
Let $E$ be an operator sequence space. Then the {\it maximal operator space\/} $\Max(E)$ over $E$ is defined by letting, for $x \in M_n(E)$,
\[
  \| x \|_{M_n(\Max(E))} := \inf \{ \| \alpha \| \| \beta \| : x = \alpha \, \diag(v_1, \ldots, v_k) \beta \},
\]
where the infimum is taken over all $k,l \in \posints$, $\alpha \in M_{n,kl}$, $\beta \in M_{k,n}$, and $v_1, \ldots, v_k$ in the closed unit ball of $E^{\fr{l}}$.
\end{definition}
\begin{remarks}
\item By \cite[Satz 4.1.10]{Lam}, $\Max(E)$ is an operator space for any operator sequence space $E$.
\item For any operator sequence space $E$, and for any operator space $F$, 
\[
  \CB(\Max(E),F) = \SB(E,C(F))
\]
holds isometrically (\cite[Satz 4.1.12]{Lam}).
\item Let $E$ be a Banach space. Then the previous remark and the isometric identity (\ref{maxsb}) combined yield that
\[
  \CB(\Max(\max(E)),F) = {\cal B}(E,M_1(F))
\]
holds isometrically for each operator space $F$, i.e.\ $\Max(\max(E)) = \MAX(E)$.
\end{remarks}
\par
There is a duality between $\Min$ and $\Max$ as between $\min$ and $\max$ and $\MIN$ and $\MAX$ (\cite[Satz 4.2.1]{Lam}):
\begin{theorem}[$\Min$-$\Max$ duality] \label{MinMax}
For any operator sequence space $E$, we have the completely isometric isomorphisms
\[
  \Min(E)^\ast = \Max(E^\ast) \qquad\text{and}\qquad \Max(E)^\ast = \Min(E^\ast)
\]
\end{theorem}
\par
We can now define the column space $\COL(E)$ and the row space $\ROW(E)$ over an {\it arbitrary\/} Banach space $E$:
\begin{definition} \label{colrowdef}
Let $E$ be a Banach space.
\begin{alphitems}
\item The {\it column space\/} over $E$ is defined as $\COL(E) := \Min(\max(E))$.
\item The {\it row space\/} over $E$ is defined as $\ROW(E) := \Max(\min(E))$.
\end{alphitems}
\end{definition}
\par
Recall that an operator space $E$ is called {\it homogeneous\/} if $\CB(E) = {\cal B}(M_1(E))$ holds isometrically.
\begin{theorem}
Let $E$ be a Banach space. Then $\COL(E)$ and $\ROW(E)$ are homogeneous operator spaces such that
\begin{equation} \label{duality}
  \COL(E)^\ast = \ROW(E^\ast) \qquad\text{and}\qquad \ROW(E)^\ast = \COL(E^\ast).
\end{equation}
\end{theorem}
\begin{proof}
Since 
\[
  {\cal B}(E) = \SB(\max(E)) = \CB(\Min(\max(E))),
\]
the homogeneity of $\COL(E)$ is clear (and similarly for $\ROW(E)$).
\par
The dualities (\ref{duality}) follow immediately from Theorems \ref{minmax} and \ref{MinMax}.
\end{proof}
\begin{remark}
It is immediate from Definition \ref{colrowdef} that, for a Banach
space $E$, 
\[
  M_{n,1}(\COL(E)) = \max(E)^{\fr{n}} = M_{n,1}(\MAX(E))
\]
and
\[
  M_{1,n}(\COL(E)) = \min(E)^{\fr{n}} = M_{1,n}(\MIN(E))
\]
holds isometrically for all $n \in \posints$. 
It follows from \cite{Math} that, for a Hilbert space $\Hilbert$, the
operator space $\COL(\Hilbert)$ in the sense of Definition
\ref{colrowdef} is the usual column Hilbert space (\cite[3.4]{ER}). An
analogous statement is true for $\ROW(\Hilbert)$.
\end{remark}
\section{Amplifying operators on $L^p$-spaces}
The following definition is from \cite{Her1}:
\begin{definition} \label{pspace}
Let $p \in (1,\infty)$. A Banach space $E$ is called a $p$-space if, for any two measure spaces $X$ and $Y$, the amplification map
\begin{equation} \label{ampl}
  {\cal B}(L^p(X),L^p(Y)) \to {\cal B}(L^p(X,E),L^p(Y,E)), \quad T \mapsto T \tensor \id_E 
\end{equation}
is an isometry.
\end{definition}
\begin{remark}
By \cite[\S4,Theorem 2]{Kwa}, a Banach space $E$ is a $p$-space if and only if it is a subspace of a quotient of an $L^p$-space. We shall, however, not require this fairly deep result, and only use Definition \ref{pspace}
and two facts from \cite{Her1}:
\begin{itemize}
\item Let $q \in [1,\infty]$. Then an $L^q$-space is a $p$-space if $p \leq q \leq 2$ or $2 \leq  q \leq p$ (\cite[Theorem 1]{Her1}). 
\item A Banach space $E$ is a $p$-space if and only if $E^\ast$ is a $p'$-space (\cite[Proposition 4]{Her1}).
\end{itemize}
\end{remark}
\par
In this section, we shall see that, for a $p$-space $E$, (\ref{ampl}) is even a complete isometry --- provided that all Banach spaces involved are equipped with their respective column space structures.
\par
We start with a proof that (\ref{ampl}) is a sequential isometry if the spaces involved are both equipped with their minimal or maximal operator sequence space structure, respectively.
\begin{lemma} \label{ampllem1}
Let $E$, $F$, and $X$ be Banach spaces. For $x \in X$, define $\pi_x^F \!: {\cal B}(E,{\cal B}(X,F)) \to {\cal B}(E,F)$ by letting
\[
  \pi_x^F(T)(y) := (Ty)(x) \qquad (T \in {\cal B}(E,{\cal B}(X,F)), \, y \in E, \, x \in X).
\]
Then the following are true:
\begin{items}
\item The equality
\[
  \| T \| = \sup \left\{ \left\| \pi_x^F(T) \right\| : x \in X, \, \| x \| \leq 1 \right\} \qquad (T \in {\cal B}(E,{\cal B}(X,F)))
\]
holds.
\item For $x \in X$ with norm one, $\pi_x$ is a quotient map.
\end{items}
\end{lemma}
\begin{proof}
We have
\begin{eqnarray*}
  \| T \| & = & \sup \{ \| T y \| : y \in E, \, \| y \| \leq 1 \} \\
  & = & \sup \{ \| (T y)(x) \| : y \in E, \, \| y \| \leq 1, \, x \in X, \, \| x \| \leq 1 \} \\
  & = & \sup \left\{ \left\| \pi_x^F(T) \right\| : x \in X, \, \| x \| \leq 1 \right\},
\end{eqnarray*}
which proves (i).
\par
To prove (ii), we define an isometric right inverse of $\pi_x$ in case $\| x \| = 1$. Fix $\phi \in X^\ast$ with $\| \phi \| = \langle x, \phi \rangle = 1$. For $T \in {\cal B}(E,F)$, define
$\tilde{T} \in {\cal B}(E,{\cal B}(X,F))$ by letting
\[
  \tilde{T}y := \phi \tensor Ty \in X^\ast \tensor^\lambda F \subset {\cal B}(X,F) \qquad (y \in E).
\]
The map ${\cal B}(E,F) \ni T \mapsto \tilde{T}$ is then the desired right inverse of $\pi_x$.
\end{proof}
\begin{corollary} \label{amplcor}
Let $p \in (1,\infty)$, let $X$ and $Y$ be measure spaces, let $E$ be a $p$-space, and let $n \in \posints$. Then the amplification map
\[
  {\cal B}\left(L^p(X),\min(L^p(Y))^{\fr{n}}\right) \to {\cal B}\left(L^p(X,E),\min(L^p(Y,E))^{\fr{n}}\right), \quad T \mapsto T \tensor \id_E 
\]
is an isometry.
\end{corollary}
\begin{proof}
Let $T \in {\cal B}(L^p(X),{\cal B}(\ell^2_n,L^p(Y)))$, and fix $\xi \in \ell^2_n$ with  $\| \xi \| \leq 1$. It follows that $\pi_\xi^{L^p(Y)} \circ T \in {\cal B}(L^p(X),L^p(Y))$. Since $E$ is a $p$-space,
we have the norm equalities
\[
  \left\| \pi_\xi^{L^p(Y)} \circ T \right\| = \left\|  \left( \pi_\xi^{L^p(Y)} \circ T \right) \tensor \id_E \right\| = \left\| \pi_\xi^{L^p(Y,E)} \circ (T \tensor \id_E) \right\|.
\]
The claim then follows from Lemma \ref{ampllem1}(i).
\end{proof}
\begin{proposition} \label{minprop}
Let $p \in (1,\infty)$, let $X$ and $Y$ be measure spaces, and let $E$ be a $p$-space. Then 
\[
  \left\| T^{\fr{m}} \right\|_{{\cal B}\left(\min(L^p(X))^{\fr{m}},\min(L^p(Y))^{\fr{mn}}\right)} = \left\|  (T \tensor \id_E)^{\fr{m}} \right\|_{{\cal B}\left(\min(L^p(X,E)^{\fr{m}}),\min(L^p(Y,E))^{\fr{mn}}\right)}
\]
holds for all $m,n \in \posints$ and for all $T \in \SB\left(\min(L^p(X)),\min(L^p(Y))^{\fr{n}}\right)$. In particular, the amplification map
\[
  \SB(\min(L^p(X)),\min(L^p(Y))) \to \SB(\min(L^p(X,E)),\min(L^p(Y,E))), \quad T \mapsto T \tensor \id_E 
\]
is a sequential isometry.
\end{proposition}
\begin{proof}
Clearly, the first assertion implies the second one.
\par
Let $m,n \in \posints$. First, note that we have for all $T \in \SB\left(\min(L^p(X)),\min(L^p(Y))^{\fr{n}}\right)$:
\begin{eqnarray*}
  \lefteqn{\left\| T^{\fr{m}} \right\|_{{\cal B}({\cal B}(\ell^2_m,L^p(X)),{\cal B}(\ell^2_m \tensor^\gamma \ell^2_n, L^p(Y)))}} & & \\
  & = & \| \id_{\ell^2_m} \tensor T \|_{{\cal B}(\ell^2_m \tensor^\lambda L^p(X),\ell^2_m \tensor^\lambda {\cal B}(\ell^2_n, L^p(Y)))}, \\
  & = & \| T \|_{{\cal B}(L^p(X),{\cal B}(\ell^2_n, L^p(Y)))}, \qquad\text{by the mapping property of $\tensor^\lambda$}, \\
  & = & \| T \tensor \id_E \|_{{\cal B}(L^p(X,E),{\cal B}(\ell^2_n, L^p(Y,E)))}, \qquad\text{by Corollary \ref{amplcor}}, \\
  & = & \| \id_{\ell^2_m} \tensor T \tensor \id_E \|_{{\cal B}(\ell^2_m \tensor^\lambda L^p(X,E),\ell^2_m \tensor^\lambda {\cal B}(\ell^2_n, L^p(Y,E)))}, \\
  &   & \qquad\text{again by the mapping property of $\tensor^\lambda$}, \\
  & = & \left\| (T \tensor \id_E)^{\fr{m}} \right\|_{{\cal B}({\cal B}(\ell^2_m,L^p(X,E)),{\cal B}(\ell^2_m \tensor^\gamma \ell^2_n, L^p(Y,E)))}.
\end{eqnarray*}
Since
\[
  \pi_\xi^{L^p(Y)}\left( T^{\fr{m}} \right) = 0 \iff \pi_\xi^{L^p(Y,E)}\left( (T \tensor \id_E)^{\fr{m}} \right) = 0
\]
for all $\xi \in \ell^2_m \tensor \ell^2_n$ and for all $T \in \SB\left(\min(L^p(X)),\min(L^p(Y))^{\fr{n}}\right)$, we conclude from Lemma \ref{ampllem1}(ii) that
\begin{equation} \label{normeq}
  \left\| \pi_\xi^{L^p(Y)}\left( T^{\fr{m}} \right) \right\| = \left\| \pi_\xi^{L^p(Y,E)}\left( (T \tensor \id_E)^{\fr{m}} \right) \right\|
\end{equation}
for all $T \in \SB\left(\min(L^p(X)),\min(L^p(Y))^{\fr{n}}\right)$ and for all $\xi \in  \ell^2_m \tensor \ell^2_n$ with $\| \xi \|_{\ell^2_m \tensor^\gamma \ell^2_n} = 1$ --- and hence for all
$\xi \in  \ell^2_m \tensor \ell^2_n$. It follows, for $T \in \SB\left(\min(L^p(X)),\min(L^p(Y))^{\fr{n}}\right)$, that
\begin{eqnarray*}
  \lefteqn{\left\| T^{\fr{m}} \right\|_{{\cal B}({\cal B}(\ell^2_m,L^p(X)),{\cal B}(\ell^2_{mn}, L^p(Y)))}} & & \\
  & = & \sup \left\{ \left\| \pi_\xi^{L^p(Y)}\left( T^{\fr{m}} \right) \right\| : \xi \in \ell^2_m \tensor \ell^2_n, \, \| \xi \|_{\ell^2_{mn}} \leq 1 \right\}, \qquad\text{by Lemma \ref{ampllem1}(i)}, \\
  & = & \sup \left\{ \left\| \pi_\xi^{L^p(Y,E)}\left( (T \tensor \id_E)^{\fr{m}} \right) \right\| : \xi \in \ell^2_m \tensor \ell^2_n, \, \| \xi \|_{\ell^2_{mn}} \leq 1 \right\}, \qquad\text{by (\ref{normeq})}, \\
  & = & \left\| (T \tensor \id_E)^{\fr{m}} \right\|_{{\cal B}({\cal B}(\ell^2_m,L^p(X,E)),{\cal B}(\ell^2_{mn}, L^p(Y,E)))}, \qquad\text{again by Lemma \ref{ampllem1}(i)}.
\end{eqnarray*}
This completes the proof.
\end{proof}
\par
Together with \cite[Proposition 4]{Her1}, Theorem \ref{adjoints}, and the $\min$-$\max$ duality, Proposition \ref{minprop} yields:
\begin{corollary} \label{maxcor}
Let $p \in (1,\infty)$, let $X$ and $Y$ be measure spaces, and let $E$ be a $p$-space. Then 
\[
  \left\| T^{\fr{m}} \right\|_{{\cal B}\left(\max(L^p(X))^{\fr{m}},\max(L^p(Y))^{\fr{mn}}\right)} = \left\|  (T \tensor \id_E)^{\fr{m}} \right\|_{{\cal B}\left(\max(L^p(X,E)^{\fr{m}}),\max(L^p(Y,E))^{\fr{mn}}\right)}
\]
holds for all $m,n \in \posints$ and for all $T \in {\cal B}\left(\max(L^p(X)),\max(L^p(Y))^{\fr{n}}\right)$. In particular, the amplification map
\[
  \SB(\max(L^p(X)),\max(L^p(Y))) \to \SB(\max(L^p(X,E)),\max(L^p(Y,E))), \quad T \mapsto T \tensor \id_E 
\]
is a sequential isometry.
\end{corollary}
\par
We can now state and prove the main result of this section:
\begin{theorem} \label{amplthm}
Let $p \in (1,\infty)$, let $X$ and $Y$ be measure spaces, and let $E$ be a $p$-space. Then the amplification map
\begin{eqnarray} 
  \CB(\COL(L^p(X)),\COL(L^p(Y))) & \to & \CB(\COL(L^p(X,E)),\COL(L^p(Y,E))), \nonumber \\
  T & \mapsto & T \tensor \id_E \label{ampl3}
\end{eqnarray}
is a complete isometry.
\end{theorem}
\begin{proof}
Let $m,n \in \posints$, and let $T \in \CB(\COL(L^p(X)),M_n(\COL(L^p(Y))))$. We can amplify $T$ to an operator $T^{\fr{m}}$ from $\max(L^p(X))^{\fr{m}}$ to
${\cal B}\left(\ell^2_n,\max(L^p(Y))^{\fr{mn}} \right)$. From Lemma \ref{ampllem1} and the first part of Corollary \ref{maxcor}, we conclude that 
\begin{eqnarray*}
  \lefteqn{\left\| T^{\fr{m}} \right\|_{{\cal B}\left(\max(L^p(X))^{\fr{m}}, {\cal B}\left(\ell^2_n,\max(L^p(Y))^{\fr{mn}} \right)\right)}} & & \\
  & = & \left\| (T \tensor \id_E)^{\fr{m}} \right\|_{{\cal B}\left(\max(L^p(X,E))^{\fr{m}}, {\cal B}\left(\ell^2_n,\max(L^p(Y,E))^{\fr{mn}} \right)\right)}.
\end{eqnarray*}
\par
An almost verbatim copy of the argument used to prove Proposition \ref{minprop} yields that
\begin{eqnarray*}
  \lefteqn{\left\| T^{(m)} \right\|_{{\cal B}\left({\cal B}\left( \ell^2_m, \max(L^p(X))^{\fr{m}}\right), {\cal B}\left(\ell^2_{mn},\max(L^p(Y))^{\fr{mn}} \right)\right)}} & & \\
  & = & \left\| (T \tensor \id_E)^{(m)} \right\|_{{\cal B}\left({\cal B}\left( \ell^2_m, \max(L^p(X))^{\fr{m}}\right), {\cal B}\left(\ell^2_{mn},\max(L^p(Y))^{\fr{mn}} \right)\right)}.
\end{eqnarray*}
Consequently, 
\[
  \| T \|_\cb = \| T \tensor \id_E \|_\cb
\]
holds, so that (\ref{ampl3}) is a complete isometry.
\end{proof}
\par
Even though we won't need it in the sequel, we note the following analog of Theorem \ref{amplthm} for row spaces: It follows immediately from the theorem due to the $\COL$-$\ROW$ duality.
\begin{corollary}
Let $p \in (1,\infty)$, let $X$ and $Y$ be measure spaces, and let $E$ be a $p$-space. Then the amplification map
\begin{eqnarray*}
  \CB(\ROW(L^p(X)),\ROW(L^p(Y))) & \to & \CB(\ROW(L^p(X,E)),\ROW(L^p(Y,E))), \\
  T & \mapsto & T \tensor \id_E 
\end{eqnarray*}
is a complete isometry.
\end{corollary}
\section{Column and row space norms on tensor products}
Let $p \in [1,\infty]$, let $X$ be a measure space, and let $E$ be a Banach space. Then the algebraic tensor product $L^p(X) \tensor E$ embeds canonically into $L^p(X,E)$. The norm of $L^p(X,E)$ restricted to
$L^p(X) \tensor E$ is a cross norm. 
\par
In this section, we want to prove an operator space analog of this fact for column and row spaces.
\begin{definition} \label{crossdef1}
Let $E$ and $F$ be operator spaces, and let $C \geq 0$. A matricial norm on $E \tensor F$ is called a {\it matricial $C$-subcross norm\/} if
\[
  \| x \tensor y \|_{M_{mn}(E \tensor F)} \leq C \| x \|_{M_m(E)} \| y \|_{M_n(F)} \qquad (m,n \in \posints, \, x \in M_m(E), \, y \in M_n(F)).
\]
If $C =1$, we simply speak of a {\it matricial subcross norm\/}.
\end{definition}
\par
There is an analog of Definition \ref{crossdef1} in the category of operator sequence spaces (\cite[Chapter 3]{Lam}):
\begin{definition} \label{crossdef2}
Let $E$ and $F$ be operator sequence spaces, and let $C \geq 0$. A sequential norm on $E \tensor F$ is called a {\it sequential $C$-subcross norm\/} if
\[
  \| x \tensor y \|_{(E \tensor F)^{\fr{mn}}} \leq C \| x \|_{E^{\fr{m}}} \| y \|_{F^{\fr{n}}} \qquad \left(m,n \in \posints, \, x \in E^{\fr{m}}, \, y \in F^{\fr{n}} \right).
\]
If $C =1$, we simply speak of a {\it sequential subcross norm\/}.
\end{definition}
\par
For the definition of Grothendieck's constant, which we will denote by $\KG$, see \cite[14.4]{DF}.
\par
Our next lemma is a consequence of \cite[26.3, Proposition 1]{DF} (see the following remark on \cite[p.\ 347]{DF}):
\begin{lemma} \label{ampllem2}
Let $p,q \in [1,\infty]$, let $X$ and $Y$ be measure spaces, and let $\Hilbert$ be a Hilbert space. Then the amplification map
\[
  {\cal B}(L^p(X),L^q(Y)) \to {\cal B}(L^p(X,\Hilbert),L^q(Y,\Hilbert))), \quad T \mapsto T \tensor \id_\Hilbert
\]
is bounded and has norm at most $\KG$. 
\end{lemma}
\begin{proposition} \label{crossprop1}
Let $p \in (1,\infty)$, let $X$ be a measure space, and let $E$ be a Banach space. Then $\min(L^p(X,E))$ yields a sequential $\KG$-subcross norm on $L^p(X) \tensor E$.
\end{proposition}
\begin{proof}
Let $m,n \in \posints$, let $x \in \min(L^p(X))^{\fr{m}}$ and let $y \in \min(E)^{\fr{n}}$. Let $S \in {\cal B}(\ell^2_m,L^p(X))$ and $T \in {\cal B}(\ell^2_n,E)$ represent $x$ and $y$, respectively. 
We have to show that
\[
  \| S \tensor T \|_{{\cal B}(\ell^2_{mn},L^p(X,E))} \leq \KG \| S \|_{{\cal B}(\ell^2_m,L^p(X))} \| T \|_{{\cal B}(\ell^2_n,E)}.
\]
By Lemma \ref{ampllem2}, we have
\[
  \| S \tensor \id_{\ell^2_n} \|_{{\cal B}(\ell^2_{mn},L^p(X,\ell^2_n))} \leq \KG \| S \|_{{\cal B}(\ell^2_m,L^p(X))},
\]
and, by \cite[7.3]{DF},
\[
  \| \id_{L^p(X)} \tensor T \|_{{\cal B}(L^p(X,\ell^2_n),L^p(X,E))} = \| T \|_{{\cal B}(\ell^2_n,E)}
\]
holds. Consequently, 
\begin{eqnarray*}
  \lefteqn{\| S \tensor T \|_{{\cal B}(\ell^2_{mn},L^p(X,E))}} \\
  & \leq & \| S \tensor \id_{\ell^2_n} \|_{{\cal B}(\ell^2_{mn},L^p(X,\ell^2_n))}
  \| \id_{L^p(X)} \tensor T \|_{{\cal B}(L^p(X,\ell^2_n),L^p(X,E))} \\
  & \leq & \KG \| S \|_{{\cal B}(\ell^2_m,L^p(X))}\| T \|_{{\cal B}(\ell^2_n,E)}
\end{eqnarray*}
holds, which completes the proof.
\end{proof}
\begin{remark}
For $p \geq 2$, we even obtain a sequential subcross norm on $L^p(X,E)$: This follows from \cite[7.2, Proprosition, and 7.3]{DF}.
\end{remark}
\begin{theorem} \label{crossthm1}
Let $p \in (1,\infty)$, let $X$ be a measure space, and let $E$ be a Banach space. Then $\ROW(L^p(X,E))$ yields a matricial $\KG$-subcross norm on $L^p(X) \tensor E$.
\end{theorem}
\begin{proof}
Let $m, n \in \posints$, and let $x \in M_m(L^p(X))$ and $y \in M_n(E)$. Let $k_x, l_x, k_y \in \posints$, and let:
\begin{itemize}
\item $\alpha_x \in M_{m,k_x l_x}$, $\beta_x \in M_{k_x,m}$, and $v_1, \ldots, v_{k_x}$ belonging to the closed unit ball of $\min(L^p(X,E))^{\fr{l_x}}$ such that
\begin{equation} \label{Max1}
  x = \alpha_x \, \diag(v_1, \ldots, v_{k_x}) \beta_x;
\end{equation} 
\item $\alpha_y \in M_{n,k_y l_y}$, $\beta_y \in M_{k_y,n}$, and $u_1, \ldots, u_{k_y}$ belonging to the closed unit ball of $\min(L^p(X,E))^{\fr{l_y}}$ such that
\begin{equation} \label{Max2}
  y = \alpha_y \, \diag(u_1, \ldots, u_{k_y}) \beta_y.
\end{equation}
\end{itemize}
It follows that $\alpha_x \tensor \alpha_y \in M_{mn,(k_x k_y) (\ell_x \ell_y)}$ and that $\beta_x \tensor \beta_y \in M_{k_x k_y,mn}$. From Proposition \ref{crossprop1}, it follows that
$\KG^{-1}(v_\mu \tensor u_\nu)$ belong to the closed unit ball of $\min(L^p(X,E))^{\fr{l_x l_y}}$ for $\mu = 1, \ldots, k_x$ and $\nu = 1, \ldots, k_y$. Consequently,
\[
  \KG^{-1} (x \tensor y) = (\alpha_x \tensor \alpha_y) \, \diag\left(\KG^{-1}(v_\mu \tensor u_\nu) :  \mu = 1, \ldots, k_x, \, \nu = 1, \ldots, k_y\right)(\beta_x \tensor \beta_y) 
\]
is a representation of $\KG^{-1} (x \tensor y)$ as in the definition of $\Max(\min(L^p(X,E)))$, i.e.\ of $\ROW(L^p(X,E))$, so that
\[
  \| \KG^{-1}(x \tensor y) \|_{M_{mn}(\ROW(L^p(X,E)))} \leq | \alpha_x \tensor \alpha_y | |\beta_x \tensor \beta_y| = | \alpha_x | |\beta_x| | \alpha_y | | \beta_y|
\]
holds. Since (\ref{Max1}) and (\ref{Max2}) are representations of $x$ and $y$ as the occur in the definition of $\Max(\min(L^p(X,E)))$, we conclude that
\[
   \| x \tensor y \|_{M_{mn}(\ROW(L^p(X,E)))} \leq \KG \| x \|_{M_m(\ROW(L^p(X)))} \| y \|_{M_n(\ROW(E))}.
\]
This yields the claim.
\end{proof}
\par
We now turn to proving the analog of Theorem \ref{crossthm1} for column spaces. 
\par
First, we need two more lemmas:
\begin{lemma} \label{maxlem1}
Let $p,q \in [1,\infty]$, let $X$ and $Y$ be a measures space, and let $\Hilbert$ be a Hilbert space. Then the amplification map
\[
  \SB(\max(L^p(X)),\max(L^q(Y))) \to \SB(\max(L^p(X,\Hilbert)),\max(L^q(Y,\Hilbert))), \quad T \mapsto  T \tensor \id_\Hilbert  
\]
is sequentially bounded with $\sb$-norm at most $\KG$.
\end{lemma}
\begin{proof}
As the proof of Corollary \ref{maxcor}, except that Lemma \ref{ampllem2} instead of Definition \ref{pspace} is invoked.
\end{proof}
\begin{lemma} \label{maxlem2}
Let $p \in (1,\infty)$, let $X$ be a measure space, and let $E$ and $F$ be Banach spaces. Then the amplification map
\[
  \SB(\max(E),\max(F)) \to \SB(\max(L^p(X,E)),\max(L^p(X,F))), \quad T \mapsto \id_{L^p(X)} \tensor T  
\]
is a sequential isometry.
\end{lemma}
\begin{proof}
As the proof of Corollary \ref{maxcor}, except that \cite[7.3]{DF} instead of Definition \ref{pspace} is invoked.
\end{proof}
\begin{theorem} \label{crossthm2}
Let $p \in (1,\infty)$, let $X$ be a measure space, and let $E$ be a Banach space. Then $\COL(L^p(X,E))$ yields a matricial $\KG$-subcross norm on $L^p(X) \tensor E$.
\end{theorem}
\begin{proof}
Let $m,n \in \posints$, and let $x \in M_m(L^p(X))$ and $y \in M_n(E)$. Let the operators $S \in {\cal B}\left(\ell^2_m, \max( L^p(X))^{\fr{m}}\right)$ and $T \in {\cal B}\left(\ell^2_n,\max(E)^{\fr{n}}\right)$ 
represent $x$ and $y$, respectively. We need to show that
\begin{equation} \label{subcross}
  \| S \tensor T \|_{{\cal B}\left(\ell^2_{mn}, \max(L^p(X,E)^{\fr{mn}} \right)} \leq \KG \| S \|_{{\cal B}\left(\ell^2_m, \max( L^p(X))^{\fr{m}}\right)} \|T \|_{{\cal B}\left(\ell^2_n,\max(E)^{\fr{n}}\right)}.
\end{equation}
First note that
\begin{eqnarray}
  \lefteqn{\KG \| S \|_{{\cal B}\left(\ell^2_m, \max( L^p(X))^{\fr{m}}\right)}} & & \nonumber \\ 
  & = & \KG \| S \|_{\SB\left(\max(\ell^2_m), \max( L^p(X))^{\fr{m}}\right)}, \qquad\text{by (\ref{maxsb})},  \nonumber \\
  & \geq & \| S \tensor \id_{\ell^2_n} \|_{\SB\left(\max(\ell^2_m(\ell^2_n)), \max( L^p(X,\ell^2_n))^{\fr{m}}\right)}, \qquad\text{by Lemma \ref{maxlem1}}, \nonumber \\
  & = & \| S \tensor \id_{\ell^2_n} \|_{{\cal B}\left(\ell^2_{mn}, \max( L^p(X,\ell^2_n))^{\fr{m}}\right)}, \qquad\text{again by (\ref{maxsb})}. \label{opnorm1}
\end{eqnarray}
On the other hand, the following holds:
\begin{eqnarray}
  \lefteqn{\| T \|_{{\cal B}\left(\ell^2_n, \max(E)^{\fr{n}}\right)}} & & \nonumber \\ 
  & = & \| T \|_{\SB\left(\max(\ell^2_n), \max(E)^{\fr{n}}\right)}  \nonumber \\
  & = & \| \id_{L^p(X)} \tensor T \|_{\SB\left(\max(L^p(X,\ell^2_n)), \max(L^p(X,E))^{\fr{n}}\right)}, \qquad\text{by Lemma \ref{maxlem2}} \nonumber \\
  & = & \left\| (\id_{L^p(X)} \tensor T)^{\fr{m}} \right\|_{{\cal B}\left(\max(L^p(X,\ell^2_n)^{\fr{m}}), \max(L^p(X,E))^{\fr{mn}}\right)}. \label{opnorm2}
\end{eqnarray}
Combined, (\ref{opnorm1}) and (\ref{opnorm2}) yield (\ref{subcross}).
\end{proof}
\section{An operator space structure for Fig\`a-Talamanca--Herz algebras}
We now use the work done in the previous sections in order to define, for $p \in (1,\infty)$ and a locally compact group $G$, a canonical operator space structure on $A_p(G)$.
\begin{definition} \label{Apopdef}
Let $G$ be a locally compact group, and let $p \in (1,\infty)$.
\begin{alphitems}
\item The canonical operator space structure on $\PM_{p'}(G)$ is the one it inherits as a closed subspace of $\CB(\COL(L^{p'}(G)))$.
\item The canonical operator space structure on $A_p(G)$ is the one it inherits as the predual of $\PM_{p'}(G)$.
\end{alphitems}
\end{definition}
\begin{remark}
From this definition, it is immediate that
\begin{equation} \label{quot}
  \ROW(L^p(G)) \Tensor \COL(L^{p'}(G)) \to A_p(G), \quad \xi \tensor \eta \mapsto \xi \ast \check{\eta}
\end{equation}
is a complete quotient map.
\end{remark}
\par
Let $G$ be a locally compact group, and let $p,q \in (1,\infty)$. Then
\[
  \lambda_{p',q'} \!: G \to {\cal B}(L^{p'}(G,L^{q'}(G))), \quad x \mapsto \lambda_{p'}(x) \tensor \lambda_{q'}(x)
\]
is a strongly continuous representation of $G$ and thus yields a representation of $L^1(G)$, which we denote likewise by $\lambda_{p',q'}$. Let the $w^\ast$-closure of $\lambda_{p',q'}(L^1(G))$ in 
${\cal B}(L^{p'}(G,L^{q'}(G)))$ be denoted by $\PM_{p',q'}(G \times G)$. Then $\PM_{p',q'}(G \times G)$ inherits a canonical operator space structure from $\CB(\COL(L^{p'}(G,L^{q'}(G)))$. Consequently, its 
predual, which we denote by $A_{p,q}(G \times G)$ has a canonical operator space structure as well. In analogy with (\ref{quot}), we have a complete quotient map
\[
  \ROW(L^p(G,L^q(G))) \Tensor\COL( L^{p'}(G,L^{q'}(G))) \to A_{p,q}(G \times G), \quad \xi \tensor \eta \mapsto \xi \ast \check{\eta}.
\]
\begin{lemma} \label{Aplem1}
Let $G$ be a locally compact group, and let $p,q \in (1,\infty)$. Then there is a canonical completely bounded map from $A_p(G) \Tensor A_q(G)$ into $A_{p,q}(G \times G)$ with $\cb$-norm at most $\KG^2$.
\end{lemma}
\begin{proof}
We have a completely isometric isomorphism 
\begin{eqnarray} 
   \lefteqn{(\ROW(L^p(G)) \Tensor \COL(L^{p'}(G))) \Tensor (\ROW(L^q(G)) \Tensor \COL(L^{q'}(G)))} & & \nonumber \\
   & \cong & (\ROW(L^p(G)) \Tensor \ROW(L^q(G))) \Tensor (\COL(L^{p'}(G)) \Tensor \COL(L^{q'}(G))) \label{shuffle}
\end{eqnarray}
by \cite[Proposition 7.1.4]{ER}. Let the left hand side of (\ref{shuffle}) be denoted by $E$, and consider the diagram
\[
  \begin{CD}
  E                     @>>> \ROW(L^p(G,L^q(G))) \Tensor \COL(L^{p'}(G),L^{q'}(G))) \\
  @VVV                       @VVV \\
  A_p(G) \Tensor A_q(G) @>>> A_{p,q}(G \times G),
  \end{CD}
\]
where the top row is the composition of (\ref{shuffle}) with the canonical completely bounded maps
\begin{equation} \label{sub1}
  \ROW(L^p(G)) \Tensor \ROW(L^q(G)) \to \ROW(L^p(G,L^q(G)))
\end{equation}
and 
\begin{equation} \label{sub2}
  \COL(L^{p'}(G)) \Tensor \COL(L^{q'}(G)) \to \COL(L^{p'}(G,L^{q'}(G))), 
\end{equation}
which exist according to Theorems \ref{crossthm1} and \ref{crossthm2} and the universal property of $\Tensor$. From Theorems \ref{crossthm1} and \ref{crossthm2}, it is also clear that both
(\ref{sub1}) and (\ref{sub2})  have $\cb$-norm at most $\KG$. 
\par
Clearly, going along the top row and down the second column is a completely bounded map with $\cb$-norm at most $\KG^2$ that factors through the kernel of the first column. 
Since the first column is a complete quotient map by \cite[Proposition 7.1.7]{ER}, we obtain the bottom row, which yields a completely bounded map whose $\cb$-norm is at most $\KG^2$ and which makes the diagram commutative.
\end{proof}
\begin{lemma} \label{Aplem2}
Let $G$ be a locally compact group, and let $p,q \in (1,\infty)$ be such that $p \leq q \leq 2$ or $2 \leq q \leq p$. Then restricting functions on $G \times G$ to the diagonal subgroup yields a
complete quotient map from $A_{p,q}(G \times G)$ onto $A_p(G)$.
\end{lemma}
\begin{proof}
Due to the particular choice of $p$ and $q$, the space $L^{q'}(G)$ is a $p'$-space by \cite[Theorem 1]{Her1}. Consequently, by Theorem \ref{amplthm}, the amplification map
\[
  \CB(\COL(L^{p'}(G))) \to \CB(\COL(L^{p'}(G,L^{q'}(G)))), \quad T \mapsto T \tensor \id_{L^{q'}(G)} 
\]
is a complete isometry. Define $W \!: L^{p'}(G,L^{q'}(G)) \to L^{p'}(G,L^{q'}(G))$ by letting
\[
  (W\xi)(x,y) := \xi(x,xy) \qquad (\xi \in L^{p'}(G,L^{q'}(G)), \, x,y \in G).
\]
Then $W$ is an invertible isometry whose inverse is given by
\[
  (W^{-1}\xi)(x,y) := \xi(x,x^{-1}y) \qquad (\xi \in L^{p'}(G,L^{q'}(G)), \, x,y \in G).
\]
Since $\COL(L^{p'}(G,L^{q'}(G)))$ is a homogeneous operator space, both $W$ and $W^{-1}$ are complete isometries on $\COL(L^{p'}(G,L^{q'}(G)))$. Consequently,
\[
  \nabla \!: \CB(\COL(L^{p'}(G))) \to \CB(\COL(L^{p'}(G,L^{q'}(G)))), \quad T \mapsto W^{-1}(T \tensor \id_{L^{q'}(G)})W 
\]
is a complete isometry. A routine calculation reveals that
\begin{equation} \label{comult}
  \nabla (\lambda_{p'}(x)) = \lambda_{p'}(x) \tensor \lambda_{q'}(x) = \lambda_{p',q'}(x,x).
\end{equation}\
It follows that $\nabla(\PM_{p'}(G)) \subset \PM_{p',q'}(G\times G)$. Moreover, $\nabla |_{\PM_{p'}(G)}$ is clearly $w^\ast$-continuous, and thus has a preadjoint $\Delta \!: A_{p,q}(G \times G) \to A_p(G)$. From (\ref{comult}),
is is immediate that $\Delta$ is the restriction to the diagonal. Finally, since $\nabla |_{\PM_{p'}(G)}$ is a complete isometry, $\Delta$ is a complete quotient map.
\end{proof}
\par
Forming the composition of the completely bounded maps in Lemmas \ref{Aplem1} and \ref{Aplem2}, we obtain the main result of this section:
\begin{theorem} \label{mainthm}
Let $G$ be a locally compact group, and let $p,q \in (1,\infty)$ be such that $p \leq q \leq 2$ or $2 \leq q \leq p$. Then pointwise multiplication induces a completely bounded map of $\cb$-norm at most $\KG^2$
from $A_p(G) \Tensor A_q(G)$ into $A_p(G)$.
\end{theorem}
\par
Letting $q = p$ in Theorem \ref{mainthm}, we obtain immediately:
\begin{corollary}
Let $G$ be a locally compact group, and let $p \in (1,\infty)$. Then $A_p(G)$ is a quantized Banach algebra.
\end{corollary}
\par
If $G$ is amenable, then $A_p(G)$ has an approximate identity bounded by one (\cite[Theorem 4.10]{Pie}). Hence, we obtain the operator space version of \cite[Theorem C]{Her1}:
\begin{corollary} \label{inclcor}
Let $G$ be an amenable, locally compact group, and let $p,q \in (1,\infty)$ be such that $p \leq q \leq 2$ or $2 \leq q \leq p$. Then $A_q(G) \subset A_p(G)$ such that the inclusion is completely bounded
with $\cb$-norm at most $\KG^2$.
\end{corollary}
\section{Operator amenability for Fig\`a-Talamanca--Herz algebras}
We conclude this paper with an extension of \cite[Theorem 3.6]{Rua1} to Fig\`a-Talamanca--Herz algebras.
\par
Let $\A$ be a quantized Banach algebra. A {\it quantized Banach $\A$-bimodule\/} is an $\A$-bimodule equipped with an operator space structure such that the module operations are completely bounded.
Let $E$ be a quantized Banach $\A$-bimodule. Then the dual space $E^\ast$ of $E$ is a quantized Banach $\A$-bimodule in a canonical fashion via
\[
  \langle x, a \cdot \phi \rangle := \langle x \cdot a, \phi \rangle \quad\text{and}\quad
  \langle x, \phi \cdot a \rangle := \langle a \cdot x, \phi \rangle
  \qquad (a \in \A, \, \phi \in E^\ast, \, x \in E).
\]
\par
A {\it derivation\/} from a quantized Banach algebra $\A$ into a quantized Banach $\A$-bimodule $E$ is a completely bounded map $D \!: \A \to E$ such that
\[
  D(ab) = a \cdot Db + (Da) \cdot b \qquad (a,b \in \A).
\]
The derivation is called {\it inner\/} if there is $x \in E$ such that
\[
  Da = a \cdot x - x \cdot a \qquad (a \in \A).
\]
\par
The following definition was introduced by Z.-J.\ Ruan in \cite{Rua1} (for completely contractive Banach algebras) and adds operator space overtones to B.\ E.\ Johnson's definition of an amenable Banach algebra 
(\cite{Joh1}):
\begin{definition}
A quantized Banach algebra $\A$ is called {\it operator amenable\/} if, for every quantized Banach $\A$-bimodule $E$, every (completely bounded) derivation $D \!: \A \to E^\ast$ is inner.
\end{definition}
\par
The classical analog of the following lemma is well known (\cite[Proposition 2.3.1]{LoA}), and the proof carries over to the quantized setting with only standard modifications. (It was formulated for completely 
contractive Banach algebras as \cite[Propositon 2.2]{Rua2}, but nowhere in the proof, complete contractivity of the multiplication is actually required.)
\begin{lemma} \label{herlem}
Let $\A$ and $\B$ be quantized Banach algebras such that $\A$ is operator amenable, and let $\theta \!: \A \to \B$ be a completely bounded algebra homomorphism with dense range. Then $\B$ is operator amenable.
\end{lemma}
\par
For Fig\`a-Talamanca--Herz algebras with their canonical operator space structure, we eventually obtain:
\begin{theorem} \label{opam}
The following are equivalent for a locally compact group $G$:
\begin{items}
\item $G$ is amenable.
\item $A(G)$ is operator amenable.
\item $A_p(G)$ is operator amenable for each $p \in (1,\infty)$.
\item There is $p \in (1,\infty)$ such that $A_p(G)$ is operator amenable.
\end{items}
\end{theorem}
\begin{proof}
(i) $\Longleftrightarrow$ (ii) is \cite[Theorem 3.6]{Rua1}.
\par
(ii) $\Longrightarrow$ (iii): Let $p \in (1,\infty)$. If $A(G)$ is operator amenable, then $G$ is amenable, so that $A(G) \subset A_p(G)$, where the inclusion is completely bounded and has 
dense range (by Corollary \ref{inclcor}). By Lemma \ref{herlem}, this yields the operator amenability of $A_p(G)$.
\par
(iii) $\Longrightarrow$ (iv) is trivial.
\par
(iv) $\Longrightarrow$ (i): Let $p \in (1,\infty)$ be such that $A_p(G)$ is operator amenable. By \cite[Proposition 2.3]{Rua1}, $A_p(G)$ then has a bounded approximate identity. This is enough
to guarantee the amenability of $G$ (\cite[Theorem 4.10]{Pie}).
\end{proof}
\begin{remarks}
\item Virtually all concepts from Banach homology can be provided with operator space overtones. For the Fourier algebra $A(G)$ of a locally compact group, this quantized Banach homology seems to be the
appropriate one when it comes to characterizing properties of $G$ in terms of cohomological properties of $A(G)$: 
\begin{itemize}
\item $A(G)$ is always operator weakly amenable (\cite{Spr}); 
\item $A(G)$ is operator biprojective if and only if $G$ is discrete (\cite{Ari}, \cite{Woo});
\item For many locally compact groups --- among them all $[\operatorname{SIN}]$-groups ---, $A(G)$ is operator biflat (\cite{RX}). The operator biflatness of $A(G)$ is systematically
investigated in \cite{ARS}.
\end{itemize}
It would be interesting to know which of these results extend to general Fig\`a-Talamanca--Herz algebras.
\item In \cite{Run3}, the third-named author showed that the Fourier algebra of a locally compact group $G$ is amenable (in the classical sense) if and only if $G$ has an abelian subgroup of finite index. The proof is mostly
operator space theoretic. It is easy to see that, if $G$ has an abelian subgroup of finite index, then $A_p(G)$ is amenable for each $p \in (1,\infty)$. In view of \cite{Run3} and Theorem \ref{opam}, it is plausible to conjecture that
$A_p(G)$ can be amenable for some $p \in (1,\infty)$ only for such $G$. It is an intriguing question, whose answer seems to be far from obvious, whether the canonical operator space structure on $A_p(G)$ --- combined
with the methods from \cite{Run3} --- can be used to affirm this conjecture.
\end{remarks}
\renewcommand{\baselinestretch}{1.0}
\dated
\vfill
\begin{tabbing}
{\it Second author's address\/}: \= Department of Mathematical and Statistical Sciences \kill 
{\it First author's address\/}: \> Fachrichtung 6.1 Mathematik \\
\> Universit\"at des Saarlandes \\
\> Postfach 151150 \\
\> 66041 Saarbr\"ucken \\
\> Germany \\[\medskipamount]
{\it E-mail\/}: \> {\tt alambert@math.uni-sb.de} \\[\medskipamount]
{\it URL\/}: \> {\tt http://www.math.uni-sb.de/$^\sim$ag-wittstock/alambert.html} \\[\bigskipamount]  
{\it Second author's address\/}: \> School of Mathematics and Statistics \\
\> 4364 Herzberg Laboratories \\ 
\> Carleton University \\ 
\> Ottawa, Ontario \\
\> Canada K1S 5B6 \\[\medskipamount]
{\it E-mail\/}: \> {\tt mneufang@math.carleton.ca} \\[\medskipamount]
{\it URL\/}: \> {\tt http://mathstat.carleton.ca/$^\sim$mneufang/} \\[\bigskipamount]  
{\it Third author's address\/}: \> Department of Mathematical and Statistical Sciences \\
\> University of Alberta \\
\> Edmonton, Alberta \\
\> Canada T6G 2G1 \\[\medskipamount]
{\it E-mail\/}: \> {\tt vrunde@ualberta.ca} \\[\medskipamount]
{\it URL\/}: \> {\tt http://www.math.ualberta.ca/$^\sim$runde/}
\end{tabbing}           
\end{document}